\documentclass[11pt,letterpaper,reqno]{amsart}

\usepackage[margin=1in]{geometry}
\usepackage{amsmath,amssymb,amsthm,mathtools}
\usepackage{hyperref}
\usepackage[nameinlink,capitalize]{cleveref}
\usepackage{enumitem, comment, amssymb}
\usepackage{color,soul}

\newtheorem{theorem}{Theorem}[section]
\newtheorem{proposition}[theorem]{Proposition}
\newtheorem{lemma}[theorem]{Lemma}

\DeclareMathOperator{\sinc}{sinc}

\DeclareMathOperator{\Id}{\mathrm{Id}}

\usepackage{scalerel}[2014/03/10]
\usepackage[usestackEOL]{stackengine}
\def\intavg{\,\ThisStyle{\ensurestackMath{%
    \stackinset{c}{0\LMpt}{c}{0\LMpt}{\SavedStyle-}{\SavedStyle\phantom{\int}}}%
    \setbox0=\hbox{$\SavedStyle\int\,$}\kern-\wd0}\int}

\theoremstyle{definition}
\newtheorem{definition}[theorem]{Definition}
\newtheorem{remark}[theorem]{Remark}

\title{Non-conical classical solutions to the one-phase free boundary problem with an isolated singularity}

\author{Minyoung Jeong}
\address{Department of Mathematics, University of California San Diego}
\email{m1jeong@ucsd.edu}

\subjclass[2020]{Primary 35R35; Secondary 49Q05, 35B40}
\keywords{Free boundary problem, one-phase problem, singularity, minimizing cone}

\begin{document}

\begin{abstract}
Several examples of one-homogeneous classical solutions to the one-phase free boundary problem with an isolated free boundary singularity are known \cite{DSJ09, H15, HKM25}. Adapting the approach of Caffarelli, Hardt, and Simon \cite{CHS84} for minimal surfaces, we construct non-conical classical solutions with an isolated free boundary singularity by perturbing a one-homogeneous solution $U$. Moreover, if $U$ is strictly stable (respectively, non-degenerate), then any solution sufficiently close to $U$ is also strictly stable (respectively, non-degenerate).
\end{abstract}

\maketitle
\markboth{}{}
\tableofcontents

\section{Introduction}\label{Introduction}
The one-phase free boundary problem is a prototypical model for stationary free interfaces arising in phase-transition and fluid-interface problems \cite{VW11, WZ12}. In this model, the interface is not prescribed a priori, but is instead determined as part of the variational problem. Mathematically, given an open set $\Omega \subseteq \mathbb R^d$ and a non-negative $g \in H^1(\Omega)$, we consider the Alt-Caffarelli energy
\begin{equation}
    \mathcal{F}(u) := \int_{\Omega} |\nabla u|^2 + \chi_{\{u > 0\}} \, dx
\end{equation}
over the admissible class
\begin{equation}
    \{u \in H^1(\Omega) : u \geq 0 \text{ a.e. in } \Omega, \quad u - g \in H^1_0(\Omega)\}.
\end{equation}

Alt and Caffarelli \cite{AC81} established the existence of minimizers by the direct method in the calculus of variations. Formally, a minimizer $u$ of this energy solves the following Euler-Lagrange equation

\begin{equation}\label{Introduction_EL_equation}
\left\{
\begin{aligned}
    \Delta u &= 0 \quad&&\text{in } \{u > 0\} \cap \Omega\\
    u &= 0 \quad&&\text{on } \partial \{u > 0\} \cap \Omega \\
    \partial_{\nu}u &= 1 \quad&&\text{on } \partial \{u > 0\} \cap \Omega\\
    u &= g \quad&&\text{on } \partial \Omega
\end{aligned}
\right.
\end{equation}
where $\nu$ denotes the inner unit normal to the free boundary $\partial \{u > 0\} \cap \Omega$.

Alt and Caffarelli \cite{AC81} also proved that minimizers of $\mathcal{F}$ are locally Lipschitz in $\Omega$ and that the free boundary $\partial\{ u > 0\} \cap \Omega$ can be decomposed into regular and singular parts, where the regular set is a $C^{1, \alpha}$ hypersurface for some $\alpha > 0$. One of their main tools is an Allard-type $\epsilon$-regularity theorem near the flat one-homogeneous minimizer: if a blow-up of $u$ at a free boundary point $x_0$ is, up to rotation, equal to $(x_d)^+$, then the free boundary is a $C^{1,\alpha}$ hypersurface in a neighborhood of $x_0$. Weiss \cite{W98, W99} subsequently proved a monotonicity formula which implies that blow-up limits at free boundary points are one-homogeneous.

Basic examples of solutions to (\ref{Introduction_EL_equation}) include half-plane and radial solutions, both of which have smooth free boundaries. Half-plane solutions are one-homogeneous, whereas radial solutions are not \cite[Section 2.3, 2.4]{V23}. Alt and Caffarelli \cite{AC81} exhibited a non-minimizing one-homogeneous solution in $\mathbb R^3$ whose free boundary has an isolated singularity. Subsequently, De Silva and Jerison \cite{DSJ09} constructed the first example of a one-homogeneous minimizer with an isolated free boundary singularity in $\mathbb R^7$. Hong \cite{H15} provided further examples of one-homogeneous solutions whose free boundaries have isolated singularities.

The examples described above with an isolated singularity are all one-homogeneous. This naturally raises the question of whether there exist solutions that are not one-homogeneous but whose free boundary have an isolated singularity.

We answer this question by perturbing a one-homogeneous classical solution $U$ with an isolated singularity at the origin. For each prescribed $\mu \geq 1$ and every sufficiently small boundary data $\varphi: \partial \{ U > 0\} \cap \mathbb S^{d-1} \rightarrow \mathbb R$, we construct a classical solution $u$ with an isolated singularity at the origin whose free boundary in $B_1$ is represented as a graph of $f : \partial \{ U > 0 \} \cap \overline{B_1} \rightarrow \mathbb R$ over the free boundary of $U$, where $f$ satisfies
\begin{equation}
    f = \varphi \text{ on $\partial \{ U > 0\} \cap \mathbb S^{d-1}$}, \quad f(x) = O(|x|^\mu) \text{ as $x \rightarrow 0$}.
\end{equation}
See Theorem \ref{Theorem_main_existence} for the precise statement. In particular, if $\mu > 1$ and $\varphi \not\equiv 0$, then $u$ is not one-homogeneous.

Furthermore, if $U$ is strictly stable (respectively, non-degenerate), then the solutions constructed above inherit strict stability (respectively, non-degeneracy), provided that $\varphi$ is sufficiently small. See Theorem \ref{Theorem_strict_stability}, \ref{Theorem_non_degeneracy}.

\subsection{Outline of the proof}
An analogous construction was done by Caffarelli, Hardt, and Simon \cite{CHS84} in minimal surface theory. Given a minimal cone $C \subseteq \mathbb R^{n+1}$ with an isolated singularity at the origin, they considered minimal surfaces represented as normal graphs
\begin{equation}
    G_u := \{x + u(x) \nu_C(x) : x \in C_1\}
\end{equation}
where we define $C_1 := (C \cap B_1) \setminus \{0\}$,  $\Gamma := C \cap \mathbb S^n$, and $u$ solves the minimal surface equation on $C_1$.

Separating the linear part of the minimal surface equation from its quadratic remainder,
\begin{equation}
    \left\{
        \begin{aligned}
            L_C u &= \mathcal{F}(u) \quad&& \text{in } C_1 \\
            u &= \varphi \quad&& \text{on } \Gamma
        \end{aligned}
    \right.    
\end{equation}
where $L_C := \Delta_C + |A_C|^2$ is the Jacobi operator on $C$, and $\mathcal{F}$ is the quadratic remainder term.

For a fixed $\mu \geq 1$, they work in weighted Holder spaces imposing
\begin{equation}
    u(x) = O(|x|^\mu)
\end{equation}
together with the corresponding decay on derivatives. For $u : \overline{C_1} \rightarrow \mathbb R$ with a sufficiently small weighted norm and solving $L_Cu = \mathcal{F}(u)$ in $C_1$, the graph $G_u$ is an embedded minimal surface with an isolated singularity at the origin. Furthermore, if $\mu > 1 $ and $u \not\equiv 0$ on $\Gamma$, then $G_u$ is non-conical.

To control the decay, they separate the radial and angular variables. We may write $u : \overline{C_1} \rightarrow \mathbb R$ as
\begin{equation}
    u(r, \omega) = \sum_{} a_j(r) \phi_j(\omega)
\end{equation}
where the eigenfunctions $\phi_j$ satisfy
\begin{equation}
    -(\Delta_\Gamma + |A_C|^2)\phi_j = \lambda_j \phi_j \quad \text{on } \Gamma.
\end{equation}
with eigenvalues, counted with multiplicity, $\lambda_1 < \lambda_2 \leq \lambda_3 \leq \cdots \rightarrow \infty$.

Then solving $L_Cu = f$ reduces to solving
\begin{equation}
    a_j''(r) + \frac{n-1}{r}a_j'(r) - \frac{\lambda_j}{r^2}a_j(r) = f_j(r), \qquad f_j(r) := \int_\Gamma f(r\omega)\phi_j(\omega)\, d\omega.
\end{equation}
The associated indicial roots are
\begin{equation}
    \gamma_j^\pm := -\frac{n-2}{2} \pm \sqrt{\left( \frac{n-2}{2} \right)^2 + \lambda_j}.
\end{equation}

Fix a $\mu > 1$ such that $\mu \ne \mathrm{Re}(\gamma_j^\pm)$ for any $j$. Let $J$ be the first integer $j$ such that $\mu < \mathrm{Re}(\gamma_j^+)$ holds. They showed that only $a_j(1)$ for $j \geq J$ can be prescribed, and $a_j(1)$ for $j < J$ are determined by the forcing term. More precisely, if we let $\Pi_J$ be the projection of $L^2(\Gamma)$ onto the closed subspace spanned by $(\phi_j)_{j \geq J}$, they solved
\begin{equation}
    \left\{
        \begin{aligned}
            Lu &= f \quad&& \text{in } C \cap B_1 \\
            \Pi_J u &= \Pi_J \varphi \quad&& \text{on } C \cap \mathbb S^n
        \end{aligned}
    \right.
\end{equation}
for each $f$ and $\varphi$.

By setting $f = \mathcal{F}(u)$, let $\mathcal{U}(u)$ be the unique solution satisfying
\begin{equation}
    \left\{
        \begin{aligned}
            L \, \mathcal{U}(u) &= \mathcal{F}(u) \quad&& \text{in } C \cap B_1 \\
            \Pi_J \mathcal{U}(u) &= \Pi_J \varphi \quad&& \text{on } C \cap \mathbb S^n
        \end{aligned}
    \right.
\end{equation}
The weighted Schauder estimate and the quadratic estimates for $\mathcal{F}$ imply that $\mathcal{U}$ is a contraction on a suitable closed ball in the weighted Holder space. Therefore, by the fixed point argument, there exists a $u : \overline{C_1} \rightarrow \mathbb R$ with $u(x) = O(|x|^\mu)$ solving
\begin{equation} \label{minimal_surface_eqaution}
    \left\{
        \begin{aligned}
            L u &= \mathcal{F}(u) \quad&& \text{in } C \cap B_1 \\
            \Pi_J u &= \Pi_J \varphi \quad&& \text{on } C \cap \mathbb S^n.
        \end{aligned}
    \right.
\end{equation}
Moreover, they showed that when $C$ is strictly stable, for sufficiently small $\varphi$, $G_u$ is also strictly stable, where $u$ solves (\ref{minimal_surface_eqaution}).

We adapt this construction to the one-phase free boundary problem. Let $U : \mathbb R^d \rightarrow [0, \infty)$ be a one-homogeneous classical solution with an isolated singularity at the origin. Define the positivity set $\Omega := \{U > 0\} \cap B_1$, and the regular free boundary $F^* := (\partial \{ U > 0\} \cap B_1) \setminus \{0\}$, and $\Sigma := \{ U > 0\} \cap \mathbb S^{d-1}$.  We seek two functions $u : B_1 \rightarrow [0, \infty)$ and $f : \partial\{U > 0\} \cap B_1 \rightarrow \mathbb R$, where $u$ solves
\begin{equation} \label{equation_01}
    \left\{
        \begin{aligned}
            \Delta u &= 0 \quad&& \text{in } \{u > 0\} \\
            \partial_{\nu_u} u &= 1 \quad&& \text{on } \partial \{u > 0\} \cap B_1
        \end{aligned}
    \right.
\end{equation}
where $\nu_u$ is the unit normal of $\partial \{u > 0\}$ pointing into $\{u > 0\}$, and the free boundary of $u$, $\partial \{u > 0\} \cap B_1$, is given by a graph of $f$ over the free boundary of $U$.

In particular, we represent the perturbed free boundary as a spherical normal graph of a function $f$ over $F^*$, and extend this perturbation to a domain deformation $\Phi_f : \Omega \rightarrow \{u > 0\}$ so that $|\Phi_f(x)| = |x|$. Pulling back (\ref{equation_01}) to $\Omega$ and writing $u\circ \Phi_f = U + w$, we reduce to a non-linear Robin problem with two unknowns $w$ and $f$ on a fixed domain $\Omega$ (Proposition \ref{Proposition_pulling_back_to_fixed_domain}). By introducing the shape derivative variable $z$, depending on $w, f$, we further reduce to a Robin problem with a single variable (Proposition \ref{Proposition_main_theorem_reduced}).

When solving the linearized problem with Robin data $h$, we first construct a bounded, linear extension $E_2(h)$ satisfying
\begin{equation}
    (\partial_\nu+ H)E_2(h) = h \quad \text{on } F^*.
\end{equation}
This extension $E_2$ maps $C_{\mu-1}^{1, \alpha}(F)$ to $C^{2, \alpha}_\mu(\Omega)$ (Lemma \ref{Lemma_extension_E2}). Subtracting this extension transfers the Robin boundary forcing to the interior. This extension only requires $h$ to be $C^{1,\alpha}$, whereas the boundary spectral construction using Steklov eigenfunction of Engelstein–Restrepo–Zhao \cite[Lemma 2.7]{ERZ25} requires $h$ to be $C^2$. For our fixed point argument, we need to solve the linearized problem for $h \in C^{1,\alpha}$.

Moreover, when proving the regularity and the estimate for our constructed solution for the linearized problem, we need a regularity result for mixed Dirichlet-Neumann boundary problem. In our case, since the two regions, one with Dirichlet boundary condition and the other with Neumann boundary condition, meet orthogonally, the regularity result can be obtained by reflection argument (Proposition \ref{Proposition_Mixed_boundary}).

A crucial distinction from the minimal surface construction is that we may prescribe any sufficiently small $\varphi : \partial \Sigma \rightarrow \mathbb R$ so that $f = \varphi$ on $\partial \Sigma$, as compared to only higher modes in the minimal surface construction. This is possible because the trace on $\Sigma$ is not prescribed. In particular, we have a linear bounded extension $E_1$ from $C^{k,\alpha}(\partial \Sigma)$ to $C^{k, \alpha}(\overline{\Sigma})$ so that $(\partial_{\nu} + H) E_1(\varphi) = 0$ on $\partial \Sigma$ and $E_1(\varphi)$ lies in the higher mode space of Robin eigenproblem (Lemma \ref{Lemma_extension_E1}).

\subsection{AI Disclosure}
AI tools were used solely for language editing and manuscript review. The mathematical content and proof strategy were developed by the author.

\subsection{Acknowledgements}
The author thanks Luca Spolaor for suggesting this problem and for his guidance and helpful discussions.

\section{Existence of non-conical solution}
\subsection{Notation and weighted Holder spaces}
All balls are centered at the origin unless otherwise specified. For $0 < r < s$, set
\begin{equation}
    A_{r, s} := B_s \setminus \overline{B}_r.
\end{equation}
For a continuous function $u : \overline{B_1} \rightarrow [0, \infty)$, define
\begin{equation}
\begin{split}
    \Omega_u &:= \{u > 0\} \cap B_1 \\
    F_u &:= \partial \Omega_u \cap B_1 \\
    F_u^* &:= F_u \setminus \{0\} \\
    \Sigma_u &:= \{ u > 0 \} \cap \mathbb S^{d-1} \\
    \partial \Sigma_u &:= \partial_{\mathbb S^{d-1}} \Sigma_u
\end{split}
\end{equation}
where $\partial \Sigma_u$ is the (topological) boundary of $\Sigma_u$ in $\mathbb S^{d-1}$.

In this article, we mainly focus on classical solutions to one-phase free boundary problem with an isolated singularity.

\begin{definition} [Classical solution  with an isolated singularity]
Let $\alpha \in (0,1)$ and $D \subseteq \mathbb R^d$ be open. We say that $u : D \rightarrow [0, \infty)$ is a classical solution to the one-phase free boundary problem in $D$ if $u \in C^{0,1}(D) \cap C^\omega(\{u>0\})$ and $\partial \{ u > 0\} \cap D$ is a (one-sided) $C^{2,\alpha}$ hypersurface, and $u$ satisfies
\begin{equation}
\begin{cases}
    \Delta u = 0 & \text{in } \{u>0\}\\
    \partial_\nu u = 1 & \text{on } \partial \{ u > 0\} \cap D
\end{cases}
\end{equation} in the classical sense, where $\nu$ denotes the inner unit normal of $\partial \{ u > 0\} \cap D$.

We say that $u$ is a classical solution to the one-phase free boundary problem with an isolated singularity at $x_0 \in D$ if $u \in C^{0,1}(D)$, its restriction to $D \setminus \{x_0\}$ is a classical solution, and $x_0 \in \partial \{ u > 0\} \cap D$ is a singular point, meaning that $\partial \{ u > 0\} \cap D$ is not a (one-sided) $C^{2,\alpha}$ hypersurface in any neighborhood of $x_0$.

When $D = \mathbb R^d$, a classical solution $u$, with or without an isolated singularity, is conical if $u$ is one-homogeneous with respect to some $a \in \mathbb R^d$, that is, $u(rx+a) = ru(x + a)$ for any $x \in \mathbb R^d$ and $r > 0$. In particular, if $u$ is conical and has an isolated singularity at $x_0$, then its center of homogeneity must be $x_0$.
\end{definition}

Let $U: \mathbb R^d \rightarrow [0, \infty)$ be a conical classical solution with an isolated singularity at the origin, and write
\begin{equation}
    \Omega := \Omega_U, \quad F := F_U, \quad F^* := F^*_U, \quad \Sigma := \Sigma_U, \quad \partial \Sigma := \partial \Sigma_U.
\end{equation}
Also, let $\nu$ be the inner unit normal of the free boundary of $U$. 

We now define the space of $C^{k,\alpha}$ functions with $r^{\mu}$ decay. 
\begin{definition} [Weighted Holder space]
Let $k$ be a non-negative integer, $\alpha \in (0, 1]$, and $E$ be either open subset or a hypersurface in $\mathbb R^d$. For $h \in C^{k, \alpha}_{\mathrm{loc}}(\overline{E}\setminus \{0\})$ and $r \in (0, \frac{1}{2}]$, define the scale-invariant $C^{k, \alpha}$ norm
\begin{equation}
\begin{split}
    \|h\|_{k, \alpha; r} &:= \sum_{j = 0}^k r^{j}\sup_{x \in E \cap A_{r, 2r}} |\nabla^j h(x)| + r^{k+\alpha} \sup_{x, y \in E \cap A_{r, 2r}} \frac{|\nabla^k h(x) - \nabla^k h(y)|}{|x-y|^\alpha}
\end{split}
\end{equation}
where $\nabla$ is the usual Euclidean connection when $E$ is an open set, and it is the induced Euclidean connection when $E$ is a hypersurface.

For $\mu \in \mathbb R$, set
\begin{equation}
    C_\mu^{k,\alpha}(E) := \{h \in C_\mathrm{loc}^{k, \alpha}(\overline{E}\setminus \{0\}) : \sup_{r \in (0, \frac{1}{2}]}r^{-\mu}\|h\|_{k, \alpha; r} < \infty \}
\end{equation}
with the norm
\begin{equation}
    \|h\|_{C_\mu^{k,\alpha}(E)} := \sup_{r \in (0, \frac{1}{2}]}r^{-\mu}\|h\|_{k, \alpha; r}.
\end{equation}
\end{definition}
One can easily verify that $(C^{k, \alpha}_{\mu}(E), \|\cdot\|_{C^{k, \alpha}_{\mu}(E)})$ is a Banach space.

\subsection{Domain perturbations and the main existence theorem}
Fix $\mu \geq 1$. Next, we consider a spherical normal coordinate near $\overline{F}\setminus \{0\}$. For $x \ne 0$ and $v \in T_x(\partial B_{|x|})$, define
\begin{equation}
    \exp_x(v) := \cos\left( \frac{|v|}{|x|} \right) x + \sinc\left( \frac{|v|}{|x|} \right) v.
\end{equation}
Notice that $\exp_x$ is the usual exponential map on the sphere of radius $|x|$.

By the tubular neighborhood theorem and the conical structure of $F$, there exists $C(U) > 0$ such that 
\[
    \mathcal{N} := \{\exp_x(t\nu_x) : |t| < C|x|, x \in \overline{F} \setminus \{0\}\}
\]
is a spherical tubular neighborhood of $\overline{F} \setminus \{0\}$ satisfying
\begin{equation}
    \overline{F} \setminus \{0\} \subseteq \mathcal{N} \subseteq \overline{B_1}.
\end{equation}
In particular, for every $x \in \mathcal{N}$, there is a unique $x_* \in \overline{F} \setminus \{0\}$ and $t \in (-C|x_*|, C|x_*|)$ such that
\begin{equation}
    x = \exp_{x_*}(t\nu_{x_*}).
\end{equation}
Define the spherical projection, normal coordinate, and normal field by
\begin{equation}
    \mathcal{P}(x) := x_*, \quad \tau(x) := t, \quad \nu^*_x := \left. \frac{d}{ds} \exp_{x_*}(s \nu(x_*)) \right|_{s = t}
\end{equation}
which are smooth in $\mathcal{N}$.

In particular,
\begin{equation}
    \exp_x(-\tau(x) \nu_x^*) = \mathcal{P}(x) \quad \text{for $x \in \mathcal{N}$}
\end{equation}
and
\[
    \mathcal{P}(x) = x, \quad \tau(x) = 0, \quad \nu^*_x = \nu_x \quad \text{for $x\in \overline{F} \setminus \{0\}$}.
\]

Finally, we consider an extension of the free boundary perturbation to a perturbation of $\Omega$. Let $\chi : \mathbb R \rightarrow [0, \infty)$ be a smooth cutoff function so that 
\begin{equation}
    \chi = 1 \text{ on $[-C/4, C/4]$}, \quad \chi = 0 \text{ on $|x| > C/2$}.    
\end{equation}
For $f \in C_\mu^{k, \alpha}(F)$, define $X_f : \overline{B_1} \rightarrow \mathbb R^d$ by
\begin{equation} \label{Perturbation}
    X_f(x) = 
    \begin{cases}
        \chi\left( \frac{\tau(x)}{|x|} \right)f(\mathcal{P}(x))\nu_x^* & x \in \mathcal{N}, \\
        0 & x \in \overline{B_1} \setminus \mathcal{N}.
    \end{cases}
\end{equation}
Clearly, $f \mapsto X_f$ is linear and $X_f \in C^{k, \alpha}_{\mathrm{loc}}(\overline{B_1} \setminus \{0\}; \mathbb R^d) \cap C^0(\overline{B_1}; \mathbb R^d)$ with
\begin{equation}
    X_f = f\nu \quad \text{on $\overline{F} \setminus \{0\}$}.
\end{equation}
Moreover, by rescaling to a fixed annulus, 
\begin{equation} \label{E1}
\begin{split}
    \|X_f\|_{2, \alpha; r} \leq C' \|f\|_{2, \alpha; r} \quad \text{for all $r \in (0, \frac{1}{2}]$}
\end{split}
\end{equation}
where $C' > 0$ depends only on $U, \alpha$, and the cutoff $\chi$.

Define the spherical normal graph over $\overline{F}\setminus \{0\}$ by
\begin{equation}
    G(f) := \{\exp_x(f(x) \nu_x) : x \in \overline{F}\setminus \{0\} \}
\end{equation}
and define the domain deformation by
\begin{equation} \label{Domain_deformation}
    \Phi_f(x) := 
    \begin{cases}
        \exp_x(X_f(x)) & x \in \overline{B_1} \setminus \{0\} \\
        0 & x = 0
    \end{cases}
\end{equation}
Then $\Phi_f \in C_{\mathrm{loc}}^{k, \alpha}(\overline{B_1}\setminus \{0\}) \cap C^0(\overline{B_1}; \mathbb R^d)$, and provided that $\mu \geq 1$ and $\|f\|_{C_\mu^{k, \alpha}(F)}$ is sufficiently small depending only on $U$ and the cutoff $\chi$, $\Phi_f$ is a $C^{k,\alpha}$ diffeomorphism on $\overline{B_1} \setminus\{0\}$.

We will assume that $\mu \geq 1$ and is not equal to any of the indicial roots of the corresponding Robin eigenproblem (see (\ref{indicial_root})).
\begin{theorem} \label{Theorem_main_existence}
For $\epsilon>0$, there exists $\delta_0(d, \alpha, U, \mu, \epsilon) > 0$ such that for every $\varphi \in C^{2, \alpha}(\partial \Sigma)$ with $\|\varphi\|_{C^{2, \alpha}(\partial \Sigma)} \leq \delta_0$, there exists a non-negative $u \in C^{0, 1}(\overline{B_1}) \cap C^{2, \alpha}_{\mathrm{loc}}(\overline{\Omega_u} \setminus \{0\})$ and $f \in C^{2,\alpha}_\mu(F)$, with $\|f\|_{C^{2,\alpha}_\mu(F)}, \|u \circ \Phi_f - U\|_{C^{2,\alpha}_\mu(\Omega)} \leq \epsilon$, satisfying
\begin{equation}
\left\{
\begin{aligned}
    \Delta u &= 0 \quad&&\text{in } \Omega_u\\
    \partial_{\nu_u}u &= 1 \quad&&\text{on } F_u^* \\
    f &= \varphi \quad&& \text{on } \partial \Sigma
\end{aligned}
\right.
\end{equation}
with
\begin{equation}
    \overline{F_u}\setminus \{0\} = G(f)
\end{equation}
where $\nu_u$ denotes the inward unit normal of $F_u^*$.
\end{theorem}

\subsection{Reduction to a nonlinear Robin problem}

Pulling back the problem to the fixed domain $\Omega$ via $\Phi_f$ and setting $v = u\circ \Phi_f$, it suffices to prove the following theorem.
\begin{proposition} \label{Proposition_pulling_back_to_fixed_domain}
For $\epsilon > 0$, there exists $\delta_0(d, \alpha, U, \mu, \epsilon) > 0$ such that for every $\varphi \in C^{2, \alpha}(\partial \Sigma)$ with $\|\varphi\|_{C^{2, \alpha}(\partial \Sigma)} \leq \delta_0$, there exists a non-negative $v \in C^{0, 1}(B_1) \cap C^{2, \alpha}_{\mathrm{loc}}(\overline{\Omega} \setminus \{0\})$ and $f \in C^{2,\alpha}_\mu(F)$, with $\|f\|_{C^{2,\alpha}_\mu(F)}$, $\|v - U\|_{C^{2,\alpha}_\mu(\Omega)} \leq \epsilon$, satisfying
\begin{equation} \label{Eq_on_Omega_1}
\left\{
\begin{aligned}
    \mathrm{div}( \Lambda_f \nabla v) &= 0 \quad&&\text{in } \Omega\\
    v &> 0 \quad&&\text{in } \Omega \\
    |(D\Phi_f)^{-T}\nabla v|^2 &= 1 \quad&& \text{on } F^*\\
    v &= 0 \quad&& \text{on } F^*\\
    f &= \varphi \quad&& \text{on } \partial \Sigma
\end{aligned}
\right.
\end{equation}
where $\Phi_f$ is defined in (\ref{Perturbation}) and (\ref{Domain_deformation}), and
\begin{equation}
    \Lambda_f := \mathrm{det} (D\Phi_f) \, (D\Phi_f)^{-1} (D\Phi_f)^{-T}.
\end{equation}
\end{proposition}

Now set $v = U + w$. We would like to write (\ref{Eq_on_Omega_1}) as 
\begin{equation}
    \mathcal{L}(w, f) = \mathcal{F}(w, f)
\end{equation}
where $\mathcal{L}$ is linear and $\mathcal{F}$ is a nonlinear remainder term satisfying the quadratic estimate
\begin{equation}
    \|\mathcal{F}(w, f)\| \leq C \|(w, f)\|^2
\end{equation}
for sufficiently small $(w, f)$ with the precise norms specified below.

For the interior equation, define
\begin{equation}
    \mathcal{L}_I(w, f) := \left. \frac{d}{dt} \, \mathrm{div}(\Lambda_{tf} \nabla (U + tw)) \right|_{t = 0}
\end{equation}
for $(w, f) \in C^{2, \alpha}_\mu(\Omega) \times C^{2, \alpha}_\mu(F)$.

Since $D\Phi_{tf} = \mathrm{Id} + tDX_f + O(t^2)$, we have
\begin{equation}
    \Lambda_{tf} = \mathrm{Id} + t (\mathrm{tr}(DX_f) \mathrm{Id} - (DX_f + DX_f^T)) + O(t^2)
\end{equation}
which yields
\begin{equation}
    \mathcal{L}_I(w, f) = \Delta w + \mathrm{div}((\mathrm{tr}(DX_f)\mathrm{Id} - (DX_f + DX_f^T)) \nabla U) 
\end{equation}
Since $\Delta U = 0$, we obtain
\begin{equation}
\begin{split}
    \mathrm{div}((\mathrm{tr}(DX_f)\mathrm{Id} - (DX_f + DX_f^T)) \nabla U)
    &= \partial_i (\partial_j X_f^j \partial_i U) - \partial_i (( \partial_j X_f^i + \partial_i X_f^j)\partial_j U) \\
    &= - (\Delta X_f^j\partial_j U + 2 \partial_i X_f^j \partial_{ij} U) \\
    &= - \Delta ( \nabla U \cdot X_f).
\end{split}
\end{equation}
Therefore,
\begin{equation}
    \mathcal{L}_I(w, f) = \Delta (w - \nabla U \cdot X_f).
\end{equation}

For the Neumann boundary equation, define
\begin{equation}
    \mathcal{L}_B(w, f) := \frac{1}{2} \left. \frac{d}{dt} \, (|(D\Phi_{tf})^{-T} \nabla (U + tw)|^2 - 1) \right|_{t = 0}
\end{equation}
for $(w, f) \in C^{2, \alpha}_\mu(\Omega) \times C^{2, \alpha}_\mu(F)$.

From $(D\Phi_{tf})^{-T} \nabla (U + tw) = \nabla U + t( \nabla w - DX_f^T \nabla U ) + O(t^2)$, $|\nabla U|^2 = 1$ and $\nabla U = \nu$, we obtain
\begin{equation}
    |(D\Phi_{tf})^{-T} \nabla (U + tw)|^2 = 1 + 2t (\partial_\nu w - \nu \cdot DX_f^T\nu) + O(t^2).
\end{equation}
Hence,
\begin{equation}
    \mathcal{L}_B(w, f) = \partial_\nu w - \nu \cdot DX_f^T \nu.
\end{equation}

Define
\begin{equation}
\begin{split}
    \mathcal{F}_I(w, f) &:= \mathcal{L}_I(w, f) - \mathrm{div} (\Lambda_f \nabla (U + w)) \\
    \mathcal{F}_B(w, f) &:= \mathcal{L}_B(w, f) - \frac{1}{2}(|D\Phi_f^{-T}\nabla (U+w)|^2 - 1)
\end{split}
\end{equation}
for $(w, f) \in C^{2, \alpha}_\mu(\Omega) \times C^{2, \alpha}_\mu(F)$.

Then (\ref{Eq_on_Omega_1}) is equivalent to
\begin{equation} \label{Eq_on_Omega_2}
\left\{
\begin{aligned}
    \mathcal{L}_I(w, f) &= \mathcal{F}_I(w, f) \quad&& \text{in } \Omega \\
    \mathcal{L}_B(w, f) &= \mathcal{F}_B(w, f) \quad&& \text{on } F^* \\
    w &= 0 \quad&& \text{on } F^* \\
    f &= \varphi \quad&& \text{on } \partial \Sigma
\end{aligned}
\right.
\end{equation}
with $v = U + w$ and $\|w\|_{C^{2,\alpha}_\mu(\Omega)}$ sufficiently small so that $v > 0$ in $\Omega$.

Now, we prove the correspondence between $(C^{2,\alpha}_\mu(\Omega) \cap \{w = 0 \text{ on } F\}) \times C^{2,\alpha}_\mu(F)$ and $C^{2, \alpha}_\mu(\Omega)$. For $(w, f) \in C^{2, \alpha}_{\mu}(\Omega) \times C^{2, \alpha}_{\mu}(F)$ and $r \in (0, \frac{1}{2}]$, define
\begin{equation}
\begin{split}
    \|(w,f)\|_{2, \alpha; r} &:= \|w\|_{2, \alpha; r} + \|f\|_{2, \alpha; r} \\
    \|(w, f)\|_{C^{2,\alpha}_\mu} &:= \|w\|_{C^{2,\alpha}_\mu(\Omega)} + \|f\|_{C^{2,\alpha}_\mu(F)}
\end{split}
\end{equation}

\begin{lemma} \label{Lemma_correspondence}
For $w \in C^{2,\alpha}_\mu(\Omega)$ with $w = 0$ on $F^*$ and $f \in C^{2,\alpha}_\mu(F)$, define
\begin{equation}
    \mathcal{E}(w, f) := w - \nabla U \cdot X_f.
\end{equation}
Then $\mathcal{E} : (C^{2,\alpha}_\mu(\Omega) \cap \{w = 0 \text{ on } F\}) \times C^{2,\alpha}_\mu(F) \rightarrow C^{2, \alpha}_\mu(\Omega)$ is a bounded linear isomorphism satisfying
\begin{equation}
    \|\mathcal{E}(w, f)\|_{2,\alpha; r} \leq C \|(w, f)\|_{2,\alpha; r}
\end{equation}
for $r \in (0, \frac{1}{2}]$ and the inverse is
\begin{equation}
    \mathcal{E}^{-1}(z) = (w_z, f_z)
\end{equation}
where
\begin{equation}
\begin{split}
    w_z &:= z + \nabla U \cdot X_{-z \mid_{F}}\\
    f_z &:= -z \mid_{F}.
\end{split}
\end{equation}
Moreover,
\begin{equation}
    \|(w_z, f_z)\|_{2,\alpha; r} \leq C \|z\|_{2,\alpha; r}
\end{equation}
for $r \in(0, \frac{1}{2}]$.
\end{lemma}
\begin{proof}
Proof is straightforward.
\end{proof}

Notice that for $z \in C^{2,\alpha}_\mu(\Omega)$
\begin{equation}
    \mathcal{L}_I(\mathcal{E}^{-1}(z)) = \Delta z \quad \text{in } \Omega.
\end{equation}

Let $H$ be the mean curvature of $F$ with the convention
\begin{equation} \label{mean_curvature}
    H = \mathrm{div}_F \nu.
\end{equation}
Since $\Delta U = 0$ up to $F^*$, we have $H = -D^2U[\nu, \nu]$. Now on $F^*$,
\begin{equation}
\begin{split}
    \mathcal{L}_B(\mathcal{E}^{-1}(z)) &= \partial_\nu w_z - \nu\cdot DX^T_{f_z}\nu \\
    &= \partial_\nu z - D^2U[\nu,\nu]z \\
    &= (\partial_\nu + H) z
.\end{split}
\end{equation}
Hence, abusing notation, for $z \in C^{2,\alpha}_\mu(\Omega)$, we write $\mathcal{F}_I(z)$ and $\mathcal{F}_B(z)$ to mean $\mathcal{F}_I(\mathcal{E}^{-1}(z))$ and $\mathcal{F}_B(\mathcal{E}^{-1}(z))$, respectively.

In the light of Lemma \ref{Lemma_correspondence}, (\ref{Eq_on_Omega_2}) is equivalent to
\begin{equation}
\left\{
\begin{aligned}
    \Delta z &= \mathcal{F}_I(z) \quad&& \text{in } \Omega \\
    (\partial_\nu + H) z &= \mathcal{F}_B(z) \quad&& \text{on } F^* \\
    z &= -\varphi \quad&& \text{on } \partial \Sigma
\end{aligned}
\right.
\end{equation}
with $w = w_z$ and $\|z\|_{C^{2,\alpha}_\mu(\Omega)}$ sufficiently small so that $\|w\|_{C^{2,\alpha}_\mu(\Omega)}$ is small, which guarantees the positivity of $v = U + w$ in $\Omega$. Note, in particular, that we no longer have the Dirichlet boundary condition on $F^*$. We summarize this as the following proposition.
\begin{proposition} \label{Proposition_main_theorem_reduced}
Let $\epsilon > 0$ and $\mu \geq 1$. There exist $\delta_0(d, \alpha, U, \mu, \epsilon) > 0$ and $C(d, \alpha, U, \mu) > 0$ such that for every $\varphi \in C^{2,\alpha}(\partial \Sigma)$ with $\|\varphi\|_{C^{2,\alpha}(\partial \Sigma)} \leq \delta_0$, there exists $z \in C^{2,\alpha}_\mu(\Omega)$ with $\|z\|_{C^{2,\alpha}_\mu(\Omega)} \leq \epsilon$ satisfying
\begin{equation}
\left\{
\begin{aligned}
    \Delta z &= \mathcal{F}_I(z) \quad&& \text{in } \Omega \\
    (\partial_\nu + H) z &= \mathcal{F}_B(z) \quad&& \text{on } F^* \\
    z &= \varphi \quad&& \text{on } \partial \Sigma
\end{aligned}
\right.
\end{equation}
with
\begin{equation}
\begin{split}
    \|z\|_{C^{2,\alpha}_\mu(\Omega)} &\leq C \|\varphi\|_{C^{2,\alpha}(\partial \Sigma)} \\
    \|z - H_K\varphi\| &\leq C \|\varphi\|_{C^{2,\alpha}(\partial \Sigma)}^2
\end{split}
\end{equation}
where $H_K$ will be defined below (\ref{Homogeneous_extension}).
\end{proposition}

We prove the quadratic estimates for $\mathcal{F}_I$ and $\mathcal{F}_B$.
\begin{lemma} \label{Lemma_quadratic_estimates}
There exists $C, \epsilon_0 > 0$ depending only on $d, \alpha, U$ such that
for any $r \in (0, \frac{1}{2}]$, $(w_1, f_1), (w_2, f_2) \in C^{2, \alpha}_{\mu}(\Omega) \times C^{2, \alpha}_{\mu}(F)$ with $\|(w_i,f_i)\|_{C^{2, \alpha}_{\mu}}< \epsilon_0$ for $i = 1,2$, the following estimates hold.
\begin{equation} \label{Lipshitz_estimate_for_remainder_term}
\begin{split}
    r^3\|\mathcal{F}_I(w_1, f_1) - \mathcal{F}_I(w_2, f_2)\|_{0, \alpha; r} &\leq C (\|(w_1, f_1)\|_{2, \alpha; r} + \|(w_2, f_2)\|_{2, \alpha; r}) \|(w_1, f_1) - (w_2, f_2)\|_{2, \alpha; r} \\
    r^2\|\mathcal{F}_B(w_1, f_1) - \mathcal{F}_B(w_2, f_2)\|_{1, \alpha; r} &\leq C (\|(w_1, f_1)\|_{2, \alpha; r} + \|(w_2, f_2)\|_{2, \alpha; r}) \|(w_1, f_1) - (w_2, f_2)\|_{2, \alpha; r}
\end{split}
\end{equation}
\end{lemma}
\begin{proof}
Define
\begin{equation}
\begin{split}
    L_f &:= (1 + \mathrm{tr} DX_f)\Id - (D X_f + DX_f^T) \\
    \mathcal{F}_{I, 1}(w, f) &:= \mathcal{L}_I(w, f) - \mathrm{div}(L_f\nabla(U+ w)) \\
    \mathcal{F}_{I, 2}(w, f) &:= \mathrm{div}(L_f\nabla(U+ w)) - \mathrm{div}(\Lambda_f\nabla(U+ w))
\end{split}
\end{equation}
Then $\mathcal{F}_I = \mathcal{F}_{I,1} + \mathcal{F}_{I,2}$, so it suffices to prove (\ref{Lipshitz_estimate_for_remainder_term}) for each term.

Fix $r \in (0, \frac{1}{2}]$, and set $Z_i := (w_i, f_i)$ for $i = 1, 2$. All pointwise estimates below hold on $\Omega \cap A_{r, 2r}$. We may write
\begin{equation}
    \mathcal{F}_{I,1}(w, f) = (DX_f + DX_f^T - \mathrm{tr} (DX_f) \mathrm{Id}) \cdot \nabla^2w + (\mathrm{div} (DX_f + DX_f^T) - \nabla \mathrm{tr}DX_f) \cdot \nabla w
\end{equation}
By linearity of $f \mapsto X_f$,
\begin{equation}
\begin{split}
    r^3|\mathcal{F}_{I,1}(w_1, f_1) - \mathcal{F}_{I,1}(w_2, f_2)| &\leq C (\|Z_1\|_{2, \alpha; r} + \|Z_2\|_{2, \alpha; r}) \|Z_1 - Z_2\|_{2, \alpha; r} \\
    r^{3+\alpha} [\mathcal{F}_{I,1}(w_1, f_1) - \mathcal{F}_{I,1}(w_2, f_2)]_{\alpha; A_{r, 2r}} &\leq C (\|Z_1\|_{2, \alpha; r} + \|Z_2\|_{2, \alpha; r}) \|Z_1 - Z_2\|_{2, \alpha; r}
\end{split}
\end{equation}
These inequalities yield the desired estimate for $\mathcal{F}_{I,1}$.

Next, let $\mathcal{R}_f := L_f - \Lambda_f$. Then
\begin{equation}
    \mathcal{F}_{I,2}(w, f) = \mathrm{div}(\mathcal{R}_f\nabla U) + \mathrm{div}(\mathcal{R}_f\nabla w)
\end{equation}
$\mathcal{R}_f$ satisfies
\begin{equation}
\begin{split}
    r^2|\mathcal{R}_f| &\leq C\|f\|_{2, \alpha; r}^2 \\
    r^2|D\mathcal{R}_f| &\leq C r^{-1 }\|f\|_{2, \alpha; r}^2 \\
    r^{3+\alpha} [D\mathcal{R}_f]_{\alpha; A_{r, 2r}} &\leq C \|f\|_{2, \alpha; r}^2
\end{split}
\end{equation}
\begin{equation}
\begin{split}
    r^2|\mathcal{R}_{f_1} - \mathcal{R}_{f_2}| &\leq C(\|f_1\|_{2, \alpha; r} + \|f_2\|_{2, \alpha; r}) \|f_1 - f_2\|_{2, \alpha; r} \\
    r^2|D(\mathcal{R}_{f_1} - \mathcal{R}_{f_2})| &\leq C r^{-1}(\|f_1\|_{2, \alpha; r} + \|f_2\|_{2, \alpha; r}) \|f_1 - f_2\|_{2, \alpha; r} \\
    r^{3+\alpha} [D(\mathcal{R}_{f_1} - \mathcal{R}_{f_2})]_{\alpha; A_{r, 2r}} &\leq C (\|f_1\|_{2, \alpha; r} + \|f_2\|_{2, \alpha; r}) \|f_1 - f_2\|_{2, \alpha; r}
\end{split}
\end{equation}
With the smallness assumptions on $w, f$, these inequalities yield the desired estimate ofr $\mathcal{F}_{I,2}$.

The same argument proves the estimate for $\mathcal{F}_B$.
\end{proof}

\subsection{Linearized problem}

We solve the linearized problem needed for the fixed point argument. We begin with some preliminary results.

Consider the Robin eigenvalue problem on $\Sigma$: for $k \geq 1$
\begin{equation} \label{Robin_eigenproblem_stability}
\left\{
\begin{aligned}
    -\Delta_\Sigma \phi_k &= \lambda_k \phi_k \quad&& \text{in } \Sigma \\
    (\partial_\nu + H) \phi_k &= 0 \quad&& \text{on } \partial \Sigma.
\end{aligned}
\right.
\end{equation}
The eigenfunctions $\phi_k$ are real analytic upto $\partial \Sigma$, and the eigenvalues, counted with multiplicity, satisfy
\begin{equation}
    \lambda_1 < \lambda_2 \leq \lambda_3 \leq \cdots, \quad \lambda_k \rightarrow \infty
\end{equation}
For each $k \geq 1$, define the indicial roots
\begin{equation} \label{indicial_root}
\begin{split}
    \gamma_{k}^\pm  &:= -\frac{d-2}{2} \pm \sqrt{\left(\frac{d-2}{2}\right)^2 + \lambda_k} \\
    \gamma_k &:= \gamma_{k}^+.
\end{split}
\end{equation}
For a fixed $\mu \geq 1$, assume that $\mu \ne \mathrm{Re}(\gamma_{k}^\pm)$ for every $k \geq 1$, and let $K$ be the first positive integer $k$ satisfying
\begin{equation}
    \mathrm{Re}(\gamma_{k}) > \mu.
\end{equation}

To solve the linearized problem, we construct two extension operators, $E_1$ and $E_2$ (Lemma \ref{Lemma_extension_E1}, \ref{Lemma_extension_E2}). The operator $E_1$ extends a function $\varphi$ on $\partial \Sigma$ to $\overline{\Sigma}$ so that 
\begin{equation}
    (\partial_\nu + H) E_1(\varphi) = 0 \quad \text{on } \partial \Sigma
\end{equation}
and $E_1(\varphi)$ is in the closed subspace of $L^2(\Sigma)$ spanned by $(\phi_k)_{k \geq K}$. The operator $E_2$ extends the Robin boundary data. Subtracting this extension transfers the Robin boundary forcing to the interior.

\begin{lemma} \label{Lemma_extension_E1}
There exists $C(d, k, \alpha, \mu, \Sigma) > 0$ and a bounded linear extension operator $E_1 : C^{k, \alpha}(\partial \Sigma) \rightarrow C^{k, \alpha}(\overline{\Sigma})$ satisfying 
\begin{enumerate}
    \item \label{E1_extension_1} $\|E_1(\varphi)\|_{C^{k, \alpha}(\overline{\Sigma})} \leq C \|\varphi\|_{C^{k, \alpha}(\partial \Sigma)}$,
    \item \label{E1_extension_2} $E_1(\varphi) = \varphi$, $\partial_\nu E_1(\varphi) = - H\varphi$ on $\partial \Sigma$,
    \item \label{E1_extension_3} $E_1(\varphi) \perp \phi_1, ..., \phi_{K-1}$ in $L^2(\Sigma)$
\end{enumerate}
for every $\varphi \in C^{k, \alpha}(\partial \Sigma)$.
\end{lemma}
\begin{proof}
By the Whitney extension theorem, there exists a linear extension map $E : C^{k, \alpha}(\partial \Sigma) \rightarrow C^{k, \alpha}(\overline{\Sigma})$ satisfying (\ref{E1_extension_1}) and (\ref{E1_extension_2}). Let $\chi \in C^\infty(\Sigma)$ be a non-negative smooth cutoff function which has a compact support in $\Sigma$. Define $a_{ij} := \langle \chi \phi_i, \phi_j \rangle_{L^2(\Sigma)}$ for $i, j = 1, ..., K-1$. Then
\begin{equation}
    a_{ij} = \delta_{ij} - \int_\Sigma (1-\chi)\phi_i\phi_j.
\end{equation}

For any $\epsilon_0> 0$, by choosing a suitable cutoff function $\chi$ depending on $\epsilon_0$ and $\phi_1, ..., \phi_{K-1}$, we may assume that 
\begin{equation}
    \left|\int_\Sigma (1-\chi)\phi_i\phi_j\right| \leq \epsilon_0
\end{equation}
for any $i, j = 1, ..., K-1$. Therefore, $(a_{ij})$ viewed as a $(K-1)\times (K-1)$ matrix is invertible. Now, let $b_i(\varphi) := \sum_{j=1}^{K-1} a^{ij} \langle E(\varphi), \phi_j\rangle$ and define
\begin{equation}
    E_1(\varphi) := E(\varphi) - \chi\sum_{j=1}^{K-1} b_j(\varphi) \phi_j.
\end{equation}
Then one can easily verify that $E_1: C^{k, \alpha}(\partial \Sigma) \rightarrow C^{k, \alpha}(\overline{\Sigma})$ is a linear extension map satisfying (\ref{E1_extension_1}), (\ref{E1_extension_2}), and (\ref{E1_extension_3}).
\end{proof}

Define the homogeneous extension $H_K : C^{2,\alpha}(\partial \Sigma) \rightarrow C_\mu^{2, \alpha}(\Omega)$ as
\begin{equation} \label{Homogeneous_extension}
    H_K(\varphi) := \sum_{k \geq K} \langle E_1(\varphi), \phi_k \rangle r^{\gamma_k} \phi_k
\end{equation}
which is a unique solution to
\begin{equation}
\left\{
\begin{aligned}
    \Delta z &= 0 \quad&& \text{in } \Omega \\
    (\partial_\nu + H) z &= 0 \quad&& \text{on } F \\
    z &= E_1(\varphi) \quad&& \text{on } \Sigma.
\end{aligned}
\right.
\end{equation}

\begin{lemma} \label{Lemma_extension_E2}
There exists $C(d, k, \alpha, \mu, U) > 0$ and a bounded extension operator $E_2 : C^{k, \alpha}_{\mu-1}(F) \rightarrow C^{k+1, \alpha}_{\mu}(\overline{\Omega})$ satisfying
\begin{enumerate}
    \item $\|E_2(h)\|_{C^{k+1, \alpha}_{\mu}(\overline{\Omega})} \leq C \|h\|_{C^{k, \alpha}_{\mu-1}(F)}$
    \item $E_2(h) = 0$, $\partial_\nu E_2(h) = h$ on $F$.
\end{enumerate}
\end{lemma}
\begin{proof}
By the Whitney extension theorem, there is a linear bounded extension $E : C^{k, \alpha}(F \cap A_{\frac{1}{4}, 1}) \rightarrow C^{k+1, \alpha}_c(\mathbb R^d)$ satisfying $E(h) = 0$ and $\partial_\nu E(h) = h$ on $F \cap A_{\frac{1}{4}, 1}$ with the estimate
\begin{equation}
    \|E(h)\|_{C^{k+1, \alpha}(\mathbb R^d)} \leq C \|h\|_{C^{k, \alpha}(F\cap A_{\frac{1}{4}, 1})}
\end{equation}
Let $r_i := 2^{-i}$, $A_i := A_{r_{i+2}, r_i}$ for $i \geq 0$. Define $h_i(x) := r_ih(r_i x)$ for $x \in F \cap A_{\frac{1}{4}, 1}$ and $E_i(x) := E(h_i)(\frac{x}{r_i})$ for $x \in \mathbb R^d$. Then $E_i = 0$ and $\partial_\nu E_i = h$ on $F \cap A_i$ with
\begin{equation}
    \|E_i\|_{k+1, \alpha; r_i} \leq C r_i (\|h\|_{k, \alpha; r_{i+2}} + \|h\|_{k, \alpha; r_{i+1}}).
\end{equation}

Now, we glue $E_i$ to construct $E_2$. Let $\eta : \mathbb R \rightarrow [0,1]$ be a non-negative smooth function satisfying $\eta(t) = 0$ for $t \leq 0$, $\eta(t) = 1$ for $t \geq 1$. Define $\lambda : \mathbb R \rightarrow [0,1]$ as
\begin{equation}
\lambda(t) :=
\begin{cases}
    \eta(t) & t \leq 1 \\
    1 - \eta(t-1) & t \geq 1
\end{cases}
\end{equation}
Since $\eta^{(j)}(t) = 0$ for any $j \geq 1$ when $t \in \{0, 1\}$, we conclude that $\lambda$ is smooth. Also, $\lambda(t) = 0$ for $t \in \mathbb R\setminus (0, 2)$ and 
\begin{equation} \label{}
    \lambda(t) + \lambda(t-1) = 1  \quad \text{for $t \in [1,2]$}.
\end{equation}

Define $\lambda_0(r) := \lambda(-\log_2(r)) + \lambda(-\log_2(r) + 1)$, $\lambda_i(r) := \lambda(-i - \log_2(r))$ for $i \geq 1$ and 
\begin{equation}
    E_2(h)(x) := \sum_{i=0}^\infty \lambda_i(|x|) E_i(x).
\end{equation}
Notice that this is locally a finite sum on $\mathbb R^d \setminus \{0\}$. Indeed, for $x \in A_{r_1, r_0}$, we have $E_2(h)(x) = E_0(x)$, and for $x \in A_{r_{i+1}, r_{i}}$ with $i \geq 1$, we have $E_2(h)(x) = \lambda_{i-1}(|x|)E_{i-1}(x) + \lambda_{i}(|x|)E_{i}(x)$. Moreover, $\lambda_{i-1}(|x|) + \lambda_{i}(|x|) = 1$ for $x \in A_{r_{i+1}, r_i}$. Therefore, it follows that $E_2(h) = 0$ and $\partial_\nu E_2(h) = h$ on $F$. Moreover, one can easily verify that $\|E_2(h)\|_{C^{k+1, \alpha}_{\mu}(\overline{\Omega})} \leq C \|h\|_{C^{k, \alpha}_{\mu-1}(F)}$.
\end{proof}

\begin{proposition} \label{Proposition_linearized_problem}
Let $\varphi \in C^{2, \alpha}(\partial \Sigma)$, $g \in C^{0, \alpha}_{\mu-2}(\Omega)$, $h \in C^{1, \alpha}_{\mu-1}(F)$. Then there exists $z \in C^{2, \alpha}_{\mu}(\Omega)$ solving
\begin{equation}
\left\{
\begin{aligned}
    \Delta z &= g  \quad&& \text{in } \Omega \\
    (\partial_\nu + H) z &= h \quad&& \text{on } F^* \\
    z &= \varphi \quad&& \text{on } \partial \Sigma
\end{aligned}
\right.
\end{equation}
which satisfies
\begin{equation} \label{linearized_solution_estimate}
\begin{split}
    \|z\|_{2, \alpha; r} &\leq C_1 r^{\mu}(\|g\|_{C^{0, \alpha}_{\mu-2}} + \|h\|_{C^{1, \alpha}_{\mu-1}} + \|\varphi\|_{C^{2,\alpha}(\partial \Sigma)}) \\
    \|z - H_J\varphi\|_{2,\alpha; r} &\leq C_1 r^{\mu}(\|g\|_{C^{0, \alpha}_{\mu-2}} + \|h\|_{C^{1, \alpha}_{\mu-1}})
\end{split}
\end{equation}
for $r \in (0, \frac{1}{2}]$ where $C_1(d, \alpha, U, \mu) >0$
\end{proposition}
\begin{proof}
Set $z_1 := E_2(h)$ and $g_1 := \Delta z_1$. Then $g_1 \in C_{\mu-2}^{0, \alpha}(\Omega)$ and 
\begin{equation}
\left\{
\begin{aligned}
    \Delta z_1 &= g_1 \quad&& \text{in } \Omega \\
    (\partial_\nu + H) z_1 &= h \quad&& \text{on } F^*.
\end{aligned}
\right.
\end{equation}

Next, let
\begin{equation}
    g_k(r) := \int_{\Sigma} (g-g_1)(r\omega) \phi_k(\omega)
\end{equation}
\begin{equation}
    z_2(r\omega) := \sum_{k \geq 1} \mathrm{Re} \left( r^{\gamma_k} \int_{b_k}^r s^{1-d-2\gamma_k} \int_0^s t^{d-1+\gamma_k} g_k(t)\, dt\, ds \right) \phi_k(\omega)
\end{equation}
where
\begin{equation}
\begin{cases}
    b_k = 0 & k \leq K-1 \\
    b_k = 1 & k \geq K.
\end{cases}
\end{equation}
Therefore, $z_2$ is a weak solution to
\begin{equation}
\left\{
\begin{aligned}
    \Delta z_2 &= g - g_1 \quad&& \text{in } \Omega \\
    (\partial_\nu + H) z_2 &= 0 \quad&& \text{on } F^*
\end{aligned}
\right.
\end{equation}
satisfying
\begin{equation}
    \|z_2(r \cdot)\|_{L^2(\Sigma)} \leq C r^\mu \|g - g_1\|_{C^{0,\alpha}_{\mu-2}(\Omega)}.
\end{equation}

Finally, set $z_3 := H_K(\varphi - (z_1 + z_2)\mid_{\partial \Sigma})$ which satisfy
\begin{equation}
\left\{
\begin{aligned}
    \Delta z_3 &= 0 &&\quad \text{in } \Omega \\
    (\partial_\nu + H) z_3 &= 0 &&\quad \text{on } F^* \\
    z_3 &= \varphi - (z_1 + z_2)\mid_{\partial \Sigma} &&\quad \text{on } \partial\Sigma.
\end{aligned}
\right.
\end{equation}

Therefore, by setting $z := z_1 + z_2 + z_3$ and the Schauder and the mixed-boundary estimate (Proposition \ref{Proposition_Mixed_boundary}) gives the desired regularity and the estimate.
\end{proof}

Define $\mathcal{S} : C^{0, \alpha}_{\mu-2}(\Omega) \times C^{1, \alpha}_{\mu-1}(F) \times C^{2,\alpha}(\partial \Sigma) \rightarrow C^{2,\alpha}_\mu(\Omega)$ so that $\mathcal{S}(g, h, \varphi)$ is the solution constructed in Proposition \ref{Proposition_linearized_problem} for each $(g, h, \varphi) \in C^{0, \alpha}_{\mu-2}(\Omega) \times C^{1, \alpha}_{\mu-1}(F) \times C^{2,\alpha}(\partial \Sigma)$. Notice that $\mathcal{S}$ is a bounded linear operator.

\subsection{Proof of the existence theorem}
\begin{proof} [Proof of Proposition \ref{Proposition_main_theorem_reduced} ]
By Lemma \ref{Lemma_correspondence}, \ref{Lemma_quadratic_estimates},
\begin{equation}
\begin{split}
    \|\mathcal{F}_I(z)\|_{0, \alpha; r} &\leq  C r^{\mu-2}\|z\|^2_{C^{2, \alpha}_\mu(\Omega)} \\
    \|\mathcal{F}_I(z_1) - \mathcal{F}_I(z_2)\|_{0, \alpha; r} &\leq C r^{\mu-2}(\|z_1\|_{C^{2, \alpha}_\mu(\Omega)} + \|z_2\|_{C^{2, \alpha}_\mu(\Omega)})\|z_1-z_2\|_{C^{2, \alpha}_\mu(\Omega)}
\end{split}
\end{equation}
\begin{equation}
\begin{split}
    \|\mathcal{F}_B(z)\|_{1, \alpha; r} &\leq  C r^{\mu-1}\|z\|^2_{C^{2, \alpha}_\mu(\Omega)} \\
    \|\mathcal{F}_B(z_1) - \mathcal{F}_B(z_2)\|_{1, \alpha; r} &\leq C r^{\mu-1}(\|z_1\|_{C^{2, \alpha}_\mu(\Omega)} + \|z_2\|_{C^{2, \alpha}_\mu(\Omega)})\|z_1-z_2\|_{C^{2, \alpha}_\mu(\Omega)}
\end{split}
\end{equation}
for $\|z\|_{C^{2,\alpha}_\mu(\Omega)}, \|z_1\|_{C^{2,\alpha}_\mu(\Omega)}, \|z_2\|_{C^{2,\alpha}_\mu(\Omega)} \leq \epsilon_0$.

For $z \in C^{2,\alpha}_\mu(\Omega)$, define
\begin{equation}
    \mathcal{U}(z) := \mathcal{S}(\mathcal{F}_I(z), \mathcal{F}_B(z), \varphi).
\end{equation}
By the above estimates and the elliptic estimate in Proposition \ref{Proposition_linearized_problem}, for $\|z\|_{C^{2,\alpha}_\mu(\Omega)}, \|z_1\|_{C^{2,\alpha}_\mu(\Omega)}, \|z_2\|_{C^{2,\alpha}_\mu(\Omega)} \leq \epsilon_0$,
\begin{equation}
\begin{split}
    \|\mathcal{U}(z)\|_{C^{2,\alpha}_\mu(\Omega)} &\leq C_2 (\|z\|_{C^{2,\alpha}_\mu(\Omega)}^2 + \|\varphi\|_{C^{2,\alpha}(\partial \Sigma)}) \\
    \|\mathcal{U}(z)- H_K(\varphi)\|_{C^{2, \alpha}_\mu(\Omega)} &\leq C_2 \|z\|^2_{C^{2, \alpha}_\mu(\Omega)} \\
    \|\mathcal{U}(z_1)- \mathcal{U}(z_2)\|_{C^{2, \alpha}_\mu(\Omega)} &\leq C_2 (\|z_1\|_{C^{2, \alpha}_\mu(\Omega)} + \|z_2\|_{C^{2, \alpha}_\mu(\Omega)})\|z_1-z_2\|_{C^{2, \alpha}_\mu(\Omega)}
\end{split}
\end{equation}
for some $C_2(d, \alpha, U, \mu)> 0$ where the last estimate follows from the interior and boundary Schauder estimates and since $\mathcal{U}(z_1)- \mathcal{U}(z_2) = \mathcal{S}(\mathcal{F}_I(z_1)-\mathcal{F}_I(z_2), \mathcal{F}_B(z_1) - \mathcal{F}_B(z_2), 0)$ with (\ref{linearized_solution_estimate}).

Let $z_0 := H_K\varphi$ and $z_{k+1} := \mathcal{U}(z_k)$ for $k \geq 0$. Fix any $\epsilon \in (0, \frac{1}{2C_2})$. Let $\delta := \min \{\frac{1}{C_1}, \frac{1}{4C_2^2}, \frac{\epsilon}{2C_2}, \frac{\epsilon_0}{2C_2}\}$. Suppose $\|\varphi\|_{C^{2,\alpha}(\partial \Sigma)} < \delta$. Then, 
\begin{equation}
\begin{split}
    \|z_k\|_{C^{2,\alpha}_\mu(\Omega)} &\leq 2C_2\delta \leq \epsilon \\
    \|z_{k+1} - z_k\|_{C^{2,\alpha}_\mu(\Omega)} &\leq C_2 \|z_0\|_{C^{2,\alpha}_\mu(\Omega)}^2 (2C_2\epsilon)^k
\end{split}
\end{equation}
for all $k \geq 0$.

Since $2C_2\epsilon < 1$, we conclude that $z_k \rightarrow z_\infty$ in $C^{2, \alpha}_\mu(\Omega)$ for some $z_\infty \in C^{2, \alpha}_\mu(\Omega)$. By continuity of $\mathcal{U}$, the limit $z_\infty$ is a fixed point of $\mathcal{U}$ with $\|z_\infty\|_{C^{2,\alpha}_\mu(\Omega)} < \epsilon$ and satisfies
\begin{equation}
\left\{
\begin{aligned}
    \Delta z_\infty &= \mathcal{F}_I(z_\infty) \quad&& \text{in } \Omega \\
    (\partial_\nu + H) z_\infty &= \mathcal{F}_B(z_\infty) \quad&& \text{on } F^* \\
    z_\infty &= \varphi \quad&& \text{on } \partial \Sigma
\end{aligned}
\right.
\end{equation}
with
\begin{equation}
\begin{split}
    \|z_\infty\|_{C^{2,\alpha}_\mu(\Omega)} &\leq C \|\varphi\|_{C^{2,\alpha}(\partial \Sigma)} \\
    \|z_\infty - H_K\varphi\|_{C^{2,\alpha}_\mu(\Omega)} &\leq C \|\varphi\|_{C^{2,\alpha}(\partial \Sigma)}^2.
\end{split}
\end{equation}

\end{proof}

\section{Stability and Non-degeneracy}
\subsection{Uniform strict stability near strictly stable conical solutions}
We begin by defining stability and strict stability for classical solutions with an isolated singularity.
\begin{definition}
Let $u : B_1 \rightarrow [0, \infty)$ be a classical solution with an isolated singularity at the origin. Let $H_u$ be the mean curvature of $F_u^*$ defined as $H_u := \mathrm{div}_{F_u^*} \nu_u$ where $\nu_u$ is the inner unit normal of $F_u^*$.

We say $u$ is stable if
\begin{equation}
    \mathcal{Q}_u(\phi) \geq 0
\end{equation}
for any non-negative $\phi \in C^{\infty}_c(B_1\setminus \{0\})$, where
\begin{equation}
    \mathcal{Q}_u(\phi) := \int_{\Omega_u} |\nabla \phi|^2 \, dx - \int_{F_u^*} H_u \phi^2 \, d\mathcal{H}^{d-1}.
\end{equation}

We say $u$ is strictly stable if there exists $c > 0$ such that
\begin{equation}
    \mathcal{Q}_u(\phi) \geq c\int_{\Omega_u} \frac{\phi^2}{|x|^2}\, dx
\end{equation}
for any non-negative $\phi \in C^{\infty}_c(B_1\setminus \{0\})$.
\end{definition}

\begin{remark}
The conical solutions $U_{d, k}$ studied by Firester–Tsiamis–Wang are strictly stable for $d \geq 7$ \cite[Theorem 1.1]{FTW26}. To our knowledge, no example is known of a stable conical classical solution with an isolated singularity that is not strictly stable. Thus, it is natural to ask whether every stable (or minimizing) cone is strictly stable. The analogous question remains open for area-minimizing hypercones with an isolated singularity \cite{C26}.
\end{remark}

Recall that stability is characterized by the second inner variation of Alt-Caffarelli energy $\mathcal{F}$. For a smooth compactly supported vector field $X \in C^{\infty}_c(B_1\setminus\{0\}; \mathbb R^d)$, we define $\Psi_t(x) := x + tX(x)$ which is a smooth diffeomorphism of $B_1$ for $|t| < \delta$. Then the second inner variation of $\mathcal{F}$ is given by
\begin{equation}
    \delta^2\mathcal{F}(u, B_1)[X] = \frac{d^2}{dt^2} \left. \mathcal{F}(u \circ \Psi^{-1}_t) \right|_{t = 0}.
\end{equation}

\begin{proposition}
$u$ is stable if and only if $\delta^2\mathcal{F}(u, B_1)[X] \geq 0$ for any $X \in C^{\infty}_c(B_1\setminus\{0\}; \mathbb R^d)$
\end{proposition}
\begin{proof}
By \cite[Proposition B.3]{K23}, we have
\begin{equation}
    \frac{1}{2} \delta^2\mathcal{F}(u, B_1)[X] = \mathcal{Q}_u(X \cdot \nabla u)
\end{equation}
for all $X \in C^{\infty}_c(B_1\setminus\{0\}; \mathbb R^d)$.

If $u$ is stable, then $\delta^2\mathcal{F}(u, B_1)[X] \geq 0$ for any $X \in C^{\infty}_c(B_1\setminus\{0\}; \mathbb R^d)$ follows immediately.

To prove the converse, let $\phi \in C_c^{\infty}(B_1 \setminus \{0\})$. Define the critical set
\begin{equation}
    Z := \{x \in \Omega_u : \nabla u(x) = 0\} \cap \mathrm{supp}(\phi).
\end{equation}
Since $\partial_{\nu} u = 1$, we have $Z \Subset \Omega_u$. Since $u$ is harmonic in $\Omega_u$, the critical set $Z$ has zero $H^1$-capacity \cite[Theorem 1.1]{HHHN99} in $\Omega_u$, and so we may take $\eta_i \in C^\infty(\Omega_u)$ such that $0 \leq \eta_i \leq 1$, $\eta_i = 0$ in a neighborhood of $Z$, $\eta_i \rightarrow 1$ in $H^1(\Omega_u)$.

Define $X_i := \eta_i \phi \frac{\nabla u}{|\nabla u|^2}$ which is well-defined in $\Omega_u$. Extending $X_i$ to $B_1$ and approximating by smooth, compactly supported vector field, we obtain $\mathcal{Q}_u (X_i \cdot \nabla u) \geq 0$ for all $i$. By taking the limit as $i \rightarrow \infty$, we conclude that $\mathcal{Q}_u(\phi) \geq 0$.
\end{proof}

We show that strict stability persists, with a uniform lower bound, under sufficiently small perturbations of a strictly stable conical solution.
\begin{theorem} \label{Theorem_strict_stability}
Suppose $U : \mathbb R^d \rightarrow [0, \infty)$ is a strictly stable conical classical solution with
\begin{equation}
    \mathcal{Q}_U(\phi) \geq c_0 \int_{\Omega} \frac{\phi^2}{|x|^2} \, dx
\end{equation}
for any non-negative $\phi \in C_c^\infty(B_1 \setminus \{0\})$.

Let $u, f$ be as in Theorem \ref{Theorem_main_existence}. There exists $\delta_1(d, \alpha, U, \mu, c_0) > 0$ such that if $\|f\|_{C^{2,\alpha}_\mu(F)} \leq \delta_1$, then $u$ is strictly stable, i.e.
\begin{equation}
    \mathcal{Q}_u(\phi) \geq \frac{c_0}{2} \int_{\Omega_u} \frac{\phi^2}{|x|^2} \, dx
\end{equation}
for any non-negative $\phi \in C_c^\infty(B_1 \setminus \{0\})$.
\end{theorem}

\begin{proof}
Let $\phi \in \phi \in C_c^\infty(B_1 \setminus \{0\})$ with $\phi \geq 0$. Define $\overline{\phi} := \phi \circ \Phi_f$, $\overline{H}_u := H_u \circ \Phi_f$. Then
\begin{equation}
    \mathcal{Q}_u(\phi) = \int_{\Omega} |(D\Phi_f)^{-T} \nabla \overline{\phi}|^2 \mathrm{det}(D\Phi_f) \, dx - \int_F \overline{H}_u \overline{\phi}^2 \mathrm{det}(D\Phi_f) |(D\Phi_f)^{-T} \nu|\, d\mathcal{H}^{d-1}.
\end{equation}
Since for any $w \in \mathbb R^d$,
\begin{equation}
    \left| \frac{|Aw|^2}{\mathrm{det}(A)} - \frac{|Bw|^2}{\mathrm{det}(B)} \right| \leq C |A-B| |w|^2
\end{equation}
for $A, B \in \mathbb R^{d\times d}$ so that $|A- \mathrm{Id}|, |B- \mathrm{Id}| \leq \delta_3$, where $\delta_3(d) > 0$, by ensuring $\|f\|_{C_\mu^{2, \alpha}(F)}$ is small so that $|(D\Phi_f)^{-T} - \mathrm{Id}| \leq C\|f\|_{C_\mu^{2, \alpha}(F)}  \leq \delta_3$, we get
\begin{equation}
    ||(D\Phi_f)^{-T} \nabla \overline{\phi}|^2 \mathrm{det}(D\Phi_f) - |\nabla \overline{\phi}|^2| \leq C \|f\|_{C_\mu^{2, \alpha}(F)} |\nabla\overline{\phi}|^2.
\end{equation}
Likewise,
\begin{equation}
    |\mathrm{det}(D\Phi_f) |(D\Phi_f)^{-T}\nu| - 1| \leq C \|f\|_{C_\mu^{2, \alpha}(F)}.
\end{equation}
Since $\overline{H}_u$ is the mean curvature of $G(f)$, we have
\begin{equation}
    |\overline{H}_u - H| \leq \frac{C\|f\|_{C_\mu^{2, \alpha}(F)}}{|x|}.
\end{equation}
Therefore, to sum up,
\begin{equation}
    |\mathcal{Q}_u(\phi) - \mathcal{Q}_U(\overline{\phi})| \leq C \|f\|_{C_\mu^{2, \alpha}(F)} \left( \int_{\Omega} |\nabla \overline{\phi}|^2\, dx + \int_F \frac{\overline{\phi}^2}{|x|}\, d\mathcal{H}^{d-1} \right).
\end{equation}
We have
\begin{equation}
\begin{split}
    \int_{\Omega} |\nabla \overline{\phi}|^2\, dx + \frac{1}{2}\int_F \frac{\overline{\phi}^2}{|x|}\, d\mathcal{H}^{d-1}
    &= \mathcal{Q}_U(\overline{\phi}) + \int_F H \overline{\phi}^2 + \frac{\overline{\phi}^2}{2|x|} \, d\mathcal{H}^{d-1} \\
    &\leq \mathcal{Q}_U(\overline{\phi}) + C\int_F \frac{\overline{\phi}^2}{|x|} \, d\mathcal{H}^{d-1}.
\end{split}
\end{equation}
By Lemma \ref{lemma_trace_inequality},
\begin{equation}
    \int_F \frac{\overline{\phi}^2}{|x|} \, d\mathcal{H}^{d-1} \leq \epsilon \int_{\Omega} |\nabla\overline{\phi}|^2 \, dx + C(\epsilon) \int_\Omega \frac{\overline{\phi}^2}{|x|^2}\, dx
\end{equation}
where $\epsilon := \frac{1}{2C}$.

Combining the previous two inequalities with the strict stability of $U$, we obtain
\begin{equation}
    \int_{\Omega} |\nabla \overline{\phi}|^2\,dx + \int_F \frac{\overline{\phi}^2}{|x|}\, d\mathcal{H}^{d-1} \leq \left( 2 + \frac{C}{c_0} \right) \mathcal{Q}_U(\overline{\phi}).
\end{equation}
Hence,
\begin{equation}
    |\mathcal{Q}_u(\phi) - \mathcal{Q}_U(\overline{\phi})| \leq C \|f\|_{C_\mu^{2, \alpha}(F)} \left( 2 + \frac{C}{c_0} \right) \mathcal{Q}_U(\overline{\phi})
\end{equation}
and for $\delta_1$ small depending only on $C(d, \alpha, U, \mu)$ and $c_0$,
\begin{equation}
    \mathcal{Q}_u(\phi) \geq \frac{3}{4} \mathcal{Q}_U(\overline{\phi}) \geq \frac{3c_0}{4} \int_{\Omega} \frac{\overline{\phi}^2}{|x|^2}\, dx.
\end{equation}
Since $|\Phi_f(x)| = |x|$ and $\mathrm{det}(D\Phi_f) \leq \frac{3}{2}$ provided that $\delta_1$ is small, we conclude that
\begin{equation}
    \mathcal{Q}_u(\phi) \geq \frac{3c_0}{4} \int_{\Omega} \frac{\overline{\phi}^2}{|x|^2}\, dx \geq \frac{c_0}{2} \int_{\Omega} \frac{\overline{\phi}^2}{|\Phi_f|^2} \mathrm{det}(D\Phi_f) \, dx = \frac{c_0}{2} \int_{\Omega_u} \frac{\phi^2}{|x|^2}\, dx.
\end{equation}
\end{proof}

\subsection{Uniform non-degeneracy near conical solutions}
Recall the definition of non-degenerate solution.
\begin{definition}
Let $u : D \rightarrow [0, \infty)$ be a classical solution. We say $u$ is non-degenerate in $D$ if there exists $c_0 > 0$ such that for any $x \in F_u$ and $B_r(x) \Subset D$,
\begin{equation}
    \sup_{B_{r}(x)} u \geq c_0 r
\end{equation}
\end{definition}

To begin with, we show that any conical classical solution with regular free boundary or an isolated singularity is non-degenerate in $\mathbb R^d$. We prove assuming $U$ has an isolated singularity at the origin, but the same argument works when the free boundary is regular everywhere.
\begin{theorem}\label{Theorem_conical_classical_with_isolated_sing_implies_non_deg}
Let $U : \mathbb R^d \rightarrow [0, \infty)$ be a conical classical solution with an isolated singularity at the origin. Then $U$ is non-degenerate in $\mathbb R^d$.
\end{theorem}
\begin{proof}
By 1-homogeneity of $U$ and compactness, there exists $\rho_0 \in (0, 1)$ such that for any $x \in F$,
\begin{equation} \label{Theorem_gradient_lower_bound}
    |\nabla U| \geq \frac{1}{2} \quad \text{on } B_{\rho_0 |x|}(x) \cap \Omega.
\end{equation}
Let $x \in F$ and $r > 0$. Consider when $|x| \geq r/2$. Then
\begin{equation}
    B_{\frac{\rho_0}{2}r}(x) \subseteq B_{\rho_0|x|}(x).
\end{equation}
By (\ref{Theorem_gradient_lower_bound}) and $U(x) = 0$, there is $x_0 \in B_{\frac{\rho_0r}{2}}(x)$ such that $U(x_0) \geq \frac{\rho_0}{4}r$. Therefore,
\begin{equation}
    \sup_{B_r(x)} U \geq U(x_0) \geq \frac{\rho_0}{4}r.
\end{equation}
Now, consider $|x| \leq r/2$. Then
\begin{equation}
\Omega \cap \left( \frac{r}{4} \mathbb S^{d-1}\right) \subseteq B_{r/2}(0) \subseteq B_r(x).
\end{equation}
Hence,
\begin{equation}
    \sup_{B_r(x)} U \geq \sup_{\Omega \cap \left( \frac{r}{4} \mathbb S^{d-1}\right)} U = \frac{\sup_{\Sigma} U}{4} r.
\end{equation}
\end{proof}

Next, we show that any classical solution close to a conical classical solution is uniformly non-degenerate.

\begin{theorem} \label{Theorem_non_degeneracy}
Suppose $U : \mathbb R^d \rightarrow [0, \infty)$ is a conical classical solution with an isolated singularity at the origin.
Let $u, f$ be as in Theorem \ref{Theorem_main_existence}. There exists $\delta_2(d, \alpha, U, \mu), c_0(U) > 0$ such that if $\|f\|_{C^{2,\alpha}_\mu(F)}, \|u \circ \Phi_f - U\|_{C^{2,\alpha}_\mu(\Omega)} \leq \delta_2$, then $u$ is non-degenerate, i.e.
\begin{equation}
    \sup_{B_r(x)} u \geq \frac{c_0}{2} r
\end{equation}
for any $B_r(x) \Subset B_1$ and $x \in F_u$.
\end{theorem}
\begin{proof}
By Theorem \ref{Theorem_conical_classical_with_isolated_sing_implies_non_deg}, we have $c_0 >0$ such that
\begin{equation}
    \sup_{B_{r}(x)} U \geq c_0 r
\end{equation}
for any $B_r(x) \Subset B_1$ and $x \in F$.

By choosing $\delta_2$ small, we may assume that $\Phi_f$ is a diffeomorphism. Let $w := u \circ \Phi_f - U$. Since $w = 0$ on $F^*$, by ensuring that $\delta_2$ is sufficiently small, we readily obtain 
\begin{equation}
    \sup_{B_r(x)} u\circ \Phi_f = \sup_{B_r(x)} (U + w) \geq \frac{2c_0}{3} r 
\end{equation}
for any $B_r(x) \Subset B_1$ and $x \in F$.

By making $\delta_2$ small, we have
\begin{equation}
    B_{\frac{3}{4}r}(\Phi_f^{-1}(x)) \subseteq \Phi_f^{-1}(B_r(x))
\end{equation}
and so
\begin{equation}
    \sup_{B_r(x)} u = \sup_{\Phi_f^{-1}(B_r(x))} u \circ \Phi_f \geq \sup_{B_{\frac{3}{4}r}(\Phi_f^{-1}(x))} u \circ \Phi_f \geq \frac{c_0}{2}r
\end{equation}
for $B_r(x) \Subset B_1$ and $x \in F_u$.
\end{proof}

\section{Appendix}
\subsection{Regularity of mixed Dirichlet-Neumann boundary problem}
We show that
\begin{proposition} \label{Proposition_Mixed_boundary}
Let $x_0 \in \partial \Sigma$ and $0 < R_1 < R_2 < 1$ with $B_{R_1} := B_{R_1}(x_0)$ and $B_{R_2} := B_{R_2}(x_0)$. Let $f \in C^{0, \alpha}(\overline{\Omega} \cap B_{R_2})$, $h_R \in C^{1, \alpha}((\overline{F}\setminus\{0\}) \cap B_{R_2})$, $h_D \in C^{2, \alpha}(\overline{\Sigma} \cap B_{R_2})$ where $h_R$ and $h_D$ satisfy the compatibility condition
\begin{equation}
    (\partial_\nu + H) h_D = h_R \quad \text{on $\partial \Sigma \cap B_{R_2}$}.
\end{equation}
Suppose $z \in H^1(\Omega \cap B_{R_2})$ is a weak solution to
\begin{equation}
\begin{cases}
    \Delta z = f  & \text{in $\Omega \cap B_{R_2}$} \\
    (\partial_\nu + H) z = h_R & \text{on $F^* \cap B_{R_2}$} \\
    z = h_D & \text{on $\Sigma \cap B_{R_2}$}
\end{cases}
\end{equation}
that is, $z = h_D$ on $\Sigma \cap B_{R_2}$ in the trace sense and
\[
    -\int_{\Omega \cap B_{R_2}} \nabla z \cdot \nabla \phi = \int_{\Omega \cap B_{R_2}} f \phi - \int_{F \cap B_{R_2}} (Hz-h_R) \phi
\]
for any $\phi \in C^\infty(\overline{\Omega \cap B_{R_2}})$ with $\phi = 0$ on $(\Omega \cap \partial B_{R_2}) \cup (\Sigma \cap B_{R_2})$.\\
Then $z \in C^{2, \alpha}(\overline{\Omega} \cap B_{R_1})$ with
\begin{equation}
    \|z\|_{C^{2,\alpha}(\overline{\Omega} \cap B_{R_1})} \leq C (\|z\|_{L^2(\Omega \cap B_{R_2})} + \|f\|_{C^{0, \alpha}(\overline{\Omega} \cap B_{R_2})} + \|h_R\|_{C^{1, \alpha}((\overline{F} \setminus\{0\}) \cap B_{R_2})} + \|h_D\|_{C^{2, \alpha}(\overline{\Sigma} \cap B_{R_2})})
\end{equation}
where $C(d, \alpha, U, x_0, R_1, R_2) >0$.
\end{proposition}
\begin{proof}
Let $x_0 \in \partial \Sigma$ and $\gamma : B_{r_0} \cap \{r, t = 0\} \rightarrow \partial \Sigma$ be a smooth local parametrization of $\partial \Sigma$ near $x_0$ so that $\gamma(0) = x_0$, where we consider $x' = (x_1, .., x_{d-2})$, $r = x_{d-1}$, $t = x_d$, and $x = (x', r, t)$. Define 
\[
    \Phi(x', r, t) := \exp_{(1-r)\gamma(x')} (t \nu((1-r)\gamma(x')))
\]
on $(x', r, t) \in B'$ where $B':= B_{r_0} \cap \{r, t \geq 0\}$ and we shrink $r_0 > 0$ if necessary. Notice that $\Phi$ is a local parametrization of $\overline{\Omega}$ near $x_0$ where $\Phi(\{r = 0\}) = \overline{\Sigma} \cap U_{x_0}$ and $\Phi(\{t = 0\}) = \overline{F} \cap U_{x_0}$ for some open neighborhood $U_{x_0}$ of $x_0$.

Let $g_{ij} := (D\Phi^{T}D\Phi)_{ij}$ and denote $G := (g_{ij})$ and its inverse by $g^{ij} := G^{-1}_{ij}$. By a direct computation, $g^{di} = \delta_{di}$ on $B' \cap \{t = 0\}$.
Let $g := \mathrm{det}(G)$ and $\hat{u} := u \circ \Phi$. We see that $\hat{z}$ is a weak solution to
\begin{equation}
\begin{cases}
    \mathrm{div}(\sqrt{g}\, G^{-1} \nabla \hat{z}) = \sqrt{g} \hat{f} & \text{in $\mathrm{int}(B')$} \\
    (\partial_t + \hat{H}) \hat{z} = \hat{h_R} & \text{on $B' \cap \{t = 0\}$} \\
    \hat{z} = \hat{h_D} & \text{on $B' \cap \{r = 0\}$}
\end{cases}
\end{equation}
that is,
\[
    -\int_{\mathrm{int}(B')} \sqrt{g} g^{ij} \partial_i \hat{z} \, \partial_j \phi = \int_{\mathrm{int}(B')} \sqrt{g} \hat{f} \phi - \int_{B' \cap \{t = 0\}} \sqrt{g}\,(\hat{H}\hat{z} - \hat{h_R}) \phi
\]
for any $\phi \in C^\infty (\overline{B'})$ with $\phi = 0$ on $(B' \cap \{r= 0\}) \cup (\overline{B'} \cap \{|x| = r_0\})$.

By the compatibility condition on $h_R$ and $h_D$ and Whitney extension theorem, there exists $\xi \in C^{2,\alpha}(\overline{B'})$ such that $(\partial_t + \hat{H})\xi = \hat{h_R}$ on $B' \cap \{t = 0\}$, $\xi = \hat{h_D}$ on $B' \cap \{r = 0\}$ satisfying
\[
    \|\xi\|_{C^{2,\alpha}(\overline{B'})} \leq C(\|\hat{h_R}\|_{C^{1,\alpha}(B' \cap \{t = 0\})} + \|\hat{h_D}\|_{C^{2,\alpha}(B' \cap \{r = 0\})} ).
\]
Therefore, by consider $\hat{z} - \xi$ and abusing notation (where we remove the hat), we may consider a weak solution to
\begin{equation}
\begin{cases}
    \mathrm{div}(\sqrt{g}\, G^{-1} \nabla z) = p & \text{in $\mathrm{int}(B')$} \\
    (\partial_t + H) z = 0 & \text{on $B' \cap \{t = 0\}$} \\
    \hat{z} = 0 & \text{on $B' \cap \{r = 0\}$}
\end{cases}
\end{equation}
that is,
\[
    -\int_{\mathrm{int}(B')} \sqrt{g} g^{ij} \partial_i z \, \partial_j \phi = \int_{\mathrm{int}(B')} p \phi - \int_{B' \cap \{t = 0\}} \sqrt{g}\, H z \phi
\]
where $p := \sqrt{g} f - \partial_j(\sqrt{g} g^{ij} \partial_i\xi)$ which is $C^{0, \alpha}(\overline{B'})$.

By Whitney extension theorem, we have nowhere vanishing $\mu \in C^{2,\alpha}(\overline{B'})$ such that $\partial_t \mu = H\mu$ on $B' \cap \{t = 0\}$. Let $w := \mu z$. Notice that $w$ is a weak solution to
\begin{equation}
\begin{cases}
    \mathrm{div}( A \nabla w) + b\cdot \nabla w + c w = q & \text{in $\mathrm{int}(B')$} \\
    \partial_t w = 0 & \text{on $B' \cap \{t = 0\}$} \\
    w = 0 & \text{on $B' \cap \{r = 0\}$}
\end{cases}
\end{equation}
where
\begin{equation}
    A := \sqrt{g}G^{-1}
\end{equation}
\begin{equation}
    b := -\frac{2}{\mu}A\nabla \mu
\end{equation}
\begin{equation}
    c := \frac{2}{\mu^2}A\nabla \mu\cdot \nabla \mu - \frac{1}{\mu} \mathrm{div}(A\nabla \mu)
\end{equation}
\begin{equation}
    q := \mu p.
\end{equation}
For $u : B' \rightarrow \mathbb R$, we write the even and odd reflections along $B' \cap \{t = 0\}$ as
\[
    u_e(x', r, t) := u(x', r, |t|)
\]
\[
    u_o(x', r, t) := \mathrm{sgn}(t) u(x', r, |t|).
\]
We define the reflections as follows.
\begin{equation}
\begin{cases}
    \tilde{a}^{ij} := a^{ij}_{e} & \text{for $1 \leq i, j \leq d-1$} \\
    \tilde{a}^{dd} := a^{dd}_{e} \\
    \tilde{a}^{di} := a^{di}_{o}, \tilde{a}^{id} := a^{id}_{o} & \text{for $1 \leq i \leq d-1$} \\
    \tilde{b}^i := b^i_e & \text{for $i \leq i \leq d-1$} \\
    \tilde{b}^d := b^d_o \\
    \tilde{c} := c_e \\
    \tilde{q} := q_e \\
    \tilde{w} := w_e.
\end{cases}
\end{equation}
Then $\tilde{w}$ is a weak solution to
\begin{equation}
\begin{cases}
    \mathrm{div}( \tilde{A} \nabla \tilde{w}) + \tilde{b}\cdot \nabla \tilde{w} + \tilde{c} \tilde{w} = \tilde{q} & \text{in $\mathrm{int}(\tilde{B})$} \\
    \tilde{w} = 0 & \text{on $\tilde{B} \cap \{r = 0\}$}
\end{cases}
\end{equation}
where $\tilde{B} := B_{r_0} \cap \{r \geq 0\}$.

Since $a^{di} = a^{id} = 0$ for $1 \leq i \leq d-1$ on $B' \cap \{t = 0 \}$, $\tilde{A}$ is Lipschitz on $\tilde{B}$. Also, $\tilde{b}, \tilde{c}, \tilde{q} \in L^\infty(\tilde{B})$. Therefore, by the boundary Calderon-Zygmund and Sobolev-Morrey embedding, we conclude that $\tilde{w} \in C^{1, \alpha}(\tilde{B} \cap B_{r_0/2}) \cap W^{2, s}(\tilde{B} \cap B_{r_0/2})$ for any $s > 1$. Hence, $k := q - b\cdot \nabla w - cw - \partial_{i} a^{ij} \partial_j w \in C^{0, \alpha}(B' \cap B_{r_0/2})$ and $w$ is a strong solution to the non-divergence equation
\begin{equation}
\begin{cases}
    a^{ij} \partial_{ij}w = k & \text{in $\mathrm{int}(B')$} \\
    \partial_t w = 0 & \text{on $B' \cap \{t = 0\}$} \\
    w = 0 & \text{on $B' \cap \{r = 0\}$}.
\end{cases}
\end{equation}
Reflecting again with $\tilde{k} := k_e$, boundary Schauder implies $w \in C^{2, \alpha}(B' \cap B_{r_0/4})$ with the desired estimate.
\end{proof}

\subsection{A weighted trace inequality on cones}
\begin{lemma} \label{lemma_trace_inequality}
For $\epsilon > 0$, there exists $C(\epsilon) > 0$ such that
\begin{equation}
    \int_{F} \frac{\phi^2}{|x|} \, d\mathcal{H}^{d-1} \leq \epsilon \int_{\Omega} |\nabla \phi|^2 \, dx + C(\epsilon) \int_{\Omega} \frac{\phi^2}{|x|^2} \, dx
\end{equation}
for every $\phi \in C_c^\infty(B_1\setminus \{0\})$.
\end{lemma}
\begin{proof}
By extending the inner unit normal $\nu$ on $F^*$, let $X : \overline{\Omega}\setminus\{0\} \rightarrow \mathbb R^d$ be a smooth vector field such that $X = -\nu$ on $F^*$. Since $F^*$ and $\Omega$ are cones, we may further impose that $X$ is 0-homogeneous, i.e. $X(rx) = X(x)$ for all $r > 0$ and $x \in \overline{\Omega}\setminus \{0\}$.

By the divergence theorem,
\begin{equation}
    \int_{F} \frac{\phi^2}{|x|} \, d\mathcal{H}^{d-1} = \int_{\Omega} \mathrm{div}\left(\frac{\phi^2}{|x|}X\right)\, dx.
\end{equation}
We have
\begin{equation}
\begin{split}
    \int_{\Omega} \mathrm{div}\left(\frac{\phi^2}{|x|}X\right)\, dx
    &= \int_\Omega \frac{\phi^2}{|x|}\mathrm{div}X + \phi^2 \nabla \left( \frac{1}{|x|} \right) \cdot X + \frac{2}{|x|} \phi\nabla \phi \cdot X\, dx\\
    &\leq C \int_\Omega \frac{\phi^2}{|x|^2} + |\nabla \phi|\frac{|\phi|}{|x|} \, dx \\
    &\leq C \int_\Omega \frac{\phi^2}{|x|^2} + \epsilon|\nabla \phi|^2 + \frac{\phi^2}{\epsilon|x|^2}\, dx \\
    &= \epsilon \int_\Omega |\nabla\phi|^2\, dx + C(\epsilon) \int_\Omega \frac{\phi^2}{|x|^2}\, dx.
\end{split}
\end{equation}

\end{proof}

\bibliographystyle{amsalpha}
\bibliography{references}

\end{document}